\documentclass[12pt]{article}
\usepackage{graphics,graphicx}
\usepackage{stmaryrd}
\usepackage{amssymb,amsmath,amsopn,amsfonts}
\usepackage[dvips]{color}
\usepackage{colordvi,multicol}
\usepackage{epsf}
\newfam\msbfam
\font\tenmsb=msbm10 \textfont\msbfam=\tenmsb \font\sevenmsb=msbm7
\scriptfont\msbfam=\sevenmsb \font\fivemsb=msbm5
\scriptscriptfont\msbfam=\fivemsb

\def\th#1{\vspace{1mm}\noindent{\bf #1}\quad}

\def\proof{\vspace{1mm}\noindent{\it Proof}\quad}

\numberwithin{equation}{section}

\def\bc{\begin{center}}
\def\ec{\end{center}}
\def\no{\noindent}
\def\hang{\hangindent\parindent}
\def\textindent#1{\indent\llap{\qquad #1\ \ \enspace}\ignorespaces}
\def\ref{\par\hang\textindent}

\begin{document}
\title{ {\bf Strong-Feller property for Navier-Stokes equations driven by space-time white noise
\thanks{Research supported in part  by NSFC (No.11671035, No.11401019)  and DFG through IRTG 1132 and CRC 701}\\} }
\author{  {\bf Rongchan Zhu}$^{\mbox{a,c}}$, {\bf Xiangchan Zhu}$^{\mbox{b,c},}$\thanks{Corresponding author}
\date{}
\thanks{E-mail address:
zhurongchan@126.com(R. C. Zhu), zhuxiangchan@126.com(X. C. Zhu)}\\\\
$^{\mbox{a}}$Department of Mathematics, Beijing Institute of Technology, Beijing 100081,  China\\
$^{\mbox{b}}$School of Science, Beijing Jiaotong University, Beijing 100044, China\\
$^{\mbox{c}}$Department of Mathematics, University of Bielefeld, D-33615 Bielefeld, Germany}

\maketitle

\noindent {\bf Abstract}

In this paper we prove strong Feller property for the Markov semigroups associated to the two or three dimensional  Navier-Stokes (N-S) equations driven by space-time white noise  using  the theory of regularity structures  introduced by Martin Hairer in \cite{Hai14}. This implies global well-posedness of 2D N-S equation driven by space-time white noise starting from every initial point in $C^\eta$ for $\eta\in (-\frac{1}{2},0)$.

\vspace{1mm}
\no{\footnotesize{\bf 2000 Mathematics Subject Classification AMS}:\hspace{2mm} 60H15, 82C28}
 \vspace{2mm}

\no{\footnotesize{\bf Keywords}:\hspace{2mm}  stochastic Navier-Stokes equation, regularity structure,  space-time white noise, renormalisation}

\section{Introduction}
In this paper, we consider the two dimensional (2D) and  three dimensional (3D) Navier-Stokes equation driven by space-time white noise:
 Recall that the Navier-Stokes equations describe the time
evolution of an incompressible fluid (see \cite{Te84}) and are given by
\begin{equation}\aligned\partial_t u+u\cdot \nabla u = &\nu \Delta u-\nabla p+\xi \\u(0)=&u_0, \quad div u = 0\endaligned\end{equation}
where $u( t,x)\in \mathbb{R}^d, d=2,3,$ denotes the value of the velocity field at time $t$ and position
$x$, $p(t,x)$ denotes the pressure, and $\xi(t,x)$ is an external force field acting on
the fluid. We will consider the case when $x\in \mathbb{T}^d, d=2,3$, the $d$-dimensional torus.
Our mathematical model for the driving force $\xi$ is a Gaussian field which is
white in time and  space.

Random Navier-Stokes equations, especially the stochastic 2D Navier-Stokes equation driven by trace-class noise, have been studied in many articles (see e.g. \cite{FG95}, \cite{HM06}, \cite{De13}, \cite{RZZ14} and the reference therein). In the two dimensional case, existence and uniqueness of  strong solutions have been obtained if the noisy forcing term is white in time and
colored in space. In the three dimensional case,  existence of martingale (=probabilistic weak) solutions, which form a Markov selection, have been constructed for  the stochastic 3D Navier-Stokes equation driven by trace-class noise in \cite{FR08}, \cite{DD03}, \cite{GRZ09}. Furthermore, the ergodicity has been obtained for every Markov selection of the martingale solutions if driven by non-degenerate trace-class noise (see \cite{FR08}).

This paper aims at proving the strong Feller property to the equation (1.1) when $\xi$ is space-time white noise.
Such a noise might not be relevant for the study of turbulence. However, in other cases, when a flow is subjected to an
external forcing with a very small time and space correlation length, a space-time white noise may be appropriate to model in this situation.

In the two dimensional case, the Navier-Stokes equation driven by space-time white noise has been studied in \cite{DD02}, where a unique global solution  starting from almost every point has been obtained by using the Gaussian invariant measure for the equation.

In the three dimensional case, we use the theory of regularity structures
 introduced by Martin Hairer in \cite{Hai14} and the paracontrolled distribution proposed by Gubinelli, Imkeller and  Perkowski in \cite{GIP13} and obtain  existence and uniqueness of  local solutions to the stochastic 3D Navier-Stokes equations driven by space-time white noise in \cite{ZZ15}.
 Recently, these two approaches have been successful in giving a meaning to a lot of ill-posed stochastic PDEs like the Kardar-Parisi-Zhang (KPZ) equation (\cite{KPZ86}, \cite{BG97}, \cite{Hai13}), the dynamical $\Phi_3^4$ model (\cite{Hai14}, \cite{CC13}) and so on.  By these theories   Markov semigroups generated by the solutions can be constructed.

Recently, there are some papers studying  the long time behavior of the Markov semigroup associated to  singular
stochastic PDEs. In \cite{RZZ16} Michael R\"{o}ckner and the authors of this paper show that the Markov semigroup for the dynamical $\Phi^4_2$ model converges to the $\Phi^4_2$
 measure
 by using asymptotic coupling method developed in \cite{HMS11}. By this and \cite{RZZ15} we can give some characterization of the $\Phi^4_2$ field. In \cite{TW16} the authors  obtain  the strong Feller property for the dynamical $P(\Phi)_2$ model and  also the
exponential ergodicity of the dynamical $\Phi^4_
2$ model.
In \cite{HM16}, strong Feller property for the Markov semigroups generated by a large class of singular
stochastic PDEs has been  obtained, which is a very useful ingredient
to establish the ergodicity of a given Markov process.

It is very natural to ask whether the solutions to the stochastic
Navier-Stokes equations driven by space-time white noise in dimensions $2$ and $3$
as constructed in \cite{DD02}, \cite{ZZ15} also satisfy the strong Feller property. Because of the presence of the Leray projection,
the results in \cite{HM16} cannot be applied in
 this case directly. In this paper we obtain strong Feller property of the solutions to the stochastic
Navier-Stokes equations driven by space-time white noise in dimensions $2$ and $3$ by using the abstract result developed in [17, Section 2].

This paper is organized as follows. In Section 2, we recall  the regularity structure theory and the framework to obtain local existence and uniqueness of  solutions
 to the 3D N-S equations  driven by space-time white noise. In Section 3, we apply the abstract result developed in [17, Section 2] to deduce the strong Feller property of the semigroup associated to the 3D N-S equation driven by space-time white noise. 2D case has been considered in Section 4. By using the strong Feller property  we obtain a unique global solution starting from every initial point in $\mathcal{C}^\eta, \eta\in (-\frac{1}{2},0)$ in Section 4. In Appendix we recall the abstract assumptions in \cite{HM16}.

\section{N-S equation by regularity structure theory}

\subsection{Preliminary on regularity structure theory}
In this subsection we recall some preliminaries for the  theory of regularity structures from \cite{Hai14} and \cite{HM15}. From this section we fixed a scaling $\mathfrak{s}=(\mathfrak{s}_0,1,...,1)$ of $\mathbb{R}^{d+1}$. We call $|\mathfrak{s}|=\mathfrak{s}_0+d$ scaling dimension.  We define the  associate  metric on $\mathbb{R}^{d+1}$  by
$$\|z-z'\|_\mathfrak{s}:=\sum_{i=0}^{d}|z_i-z'_i|^{1/\mathfrak{s}_i}.$$ For $k=(k_0,...,k_d)$ we define $|k|_\mathfrak{s}=\sum_{i=0}^d\mathfrak{s}_ik_i$.
\vskip.10in
\th{Definition 2.1} A regularity structure $\mathfrak{T}=(A,T,G)$ consists of the following elements:

(i) An index set $A\subset \mathbb{R}$ such that $0\in A$, $A$ is bounded from below and  locally finite.

(ii) A model space $T$, which is a graded vector space $T=\oplus_{\alpha\in A}T_\alpha$, with each $T_\alpha$ a Banach space. Furthermore, $T_0$ is one-dimensional and has a basis vector $\mathbf{1}$. Given $\tau\in T$ we write $\|\tau\|_\alpha$ for the norm of its component in $T_\alpha$.

(iii) A structure group $G$ of (continuous) linear operators acting on $T$ such that for every $\Gamma\in G$, every $\alpha\in A$ and every $\tau_\alpha\in T_\alpha$ one has
$$\Gamma\tau_\alpha-\tau_\alpha\in T_{<\alpha}:=\bigoplus_{\beta<\alpha}T_\beta.$$
Furthermore, $\Gamma\mathbf{1}=\mathbf{1}$ for every $\Gamma\in G$.
\vskip.10in
Now we have the the usual
polynomial regularity structure $\bar{T} =\bigoplus_{n\in\mathbb{N}}\bar{T}_n$ given by all polynomials in $d+1$ indeterminates, let us call them $X_0,..., X_d$, which denote the time and space directions respectively.  Denote $X^k=X_0^{k_0}\cdot\cdot\cdot X_d^{k_d}$ with $k$ a multi-index. In this case,   $A = \mathbb{N}$ and $\bar{T}_n$ denote the
space of monomials that are homogeneous of degree $n$.  The structure group can be defined by $\Gamma_hX^k=(X-h)^k$, $h\in\mathbb{R}^{d+1}$.

Given a smooth compactly supported test function $\varphi$, $x, y\in\mathbb{R}^d$, $\lambda>0$, we define
$$\varphi_x^\lambda(y)=\lambda^{-d}\varphi(\frac{y-x}{\lambda}).$$

Denote by $\mathcal{B}_r$ the set of smooth test functions $\varphi:\mathbb{R}^{d}\mapsto\mathbb{R}$ that are supported in the centred  ball of radius $1$ and such that their derivatives of order up to $r$ are uniformly bounded by $1$. We denote by $\mathcal{S}'$ the space of all distributions on $\mathbb{R}^{d}$.

 Now we give the definition of a model, which is a concrete way of associating every element in  the abstract regularity structure to the actual Taylor polynomial at every point.
For the Navier-Stokes equation we need to consider heat kernel composed with the Leray projection, which is not smooth on $\mathbb{R}^{d+1}\backslash\{0\}$. So we cannot apply [14, Lemma 5.5] directly. Instead we use  the inhomogeneous modelled distribution introduced in \cite{HM15}.
\vskip.10in

\th{Definition 2.2} Given a regularity structure $\mathfrak{T}$, an inhomogeneous model $(\Pi,\Gamma, \Sigma)$ consists of the following three elements:
\begin{itemize}
  \item  A collection of maps $\Gamma^t:\mathbb{R}^d\times \mathbb{R}^d\rightarrow \mathcal{G}$ parametrized by $t\in\mathbb{R}$, such that
  $$\Gamma^t_{xx}=1,\quad \Gamma^t_{xy}\Gamma^t_{yz}=\Gamma^t_{xz},$$
  for any $x,y,z\in\mathbb{R}^d$ and $t\in \mathbb{R}$, and the action of $\Gamma_{xy}^t$ on polynomials is given as above with $h=(0,y-x)$.
  \item A collection of maps $\Sigma_x:\mathbb{R}\times \mathbb{R}\rightarrow\mathcal{G}$, parametrized by $x\in\mathbb{R}^d$, such that, for any $x\in\mathbb{R}^d$ and $s,r,t\in\mathbb{R}$, one has
      $$\Sigma_x^{tt}=1, \quad \Sigma_x^{sr}\Sigma_x^{rt}=\Sigma_x^{st}, \quad \Sigma_x^{st}\Gamma^t_{xy}=\Gamma^s_{xy}\Sigma^{st}_y,$$
      and the action of $\Sigma_x^{st}$ on polynomials is given as above with $h=(t-s,0)$.
  \item A collection of linear maps $\Pi^t_x:T\rightarrow \mathcal{S}'$, such that
  $$\Pi^t_y=\Pi^t_x\Gamma^t_{xy}, \quad (\Pi^t_xX^{(0,\bar{k})})(y)=(y-x)^{\bar{k}},\quad  (\Pi^t_xX^{(k_0,\bar{k})})(y)=0,$$
  for all $x,y\in\mathbb{R}^d, t\in\mathbb{R}, \bar{k}\in\mathbb{N}^d, k_0\in\mathbb{N}$ such that $k_0>0$.
\end{itemize}
Moreover, for any $\gamma>0$ and every $T>0$, there is a constant $C$ for which the analytic bounds
$$|\langle \Pi^t_x\tau, \varphi^\lambda_x\rangle|\leq C\|\tau\|_l\lambda^l,\quad \|\Gamma_{xy}^t\tau\|_m\leq C\|\tau\|_l|x-y|^{l-m},  $$
$$\|\Sigma_x^{st}\tau\|_m\leq C\|\tau\|_l|t-s|^{(l-m)/\mathfrak{s}_0},$$
holds uniformly over all $\tau\in T_l$,  $l\in A$ with $l<\gamma$, all $m\in A$ such that $m<l$, and all test functions $\varphi\in \mathcal{B}_r$ with $r>-\inf A$, and all $t,s\in[-T,T]$ and $x,y\in\mathbb{R}^d$ such that $|t-s|\leq 1$ and $|x-y|\leq 1$.

\vskip.10in
For a model $Z=(\Pi,\Gamma,\Sigma)$ we denote by $\|\Pi\|_{\gamma;T}, \|\Gamma\|_{\gamma;T}$ and $\|\Sigma\|_{\gamma;T}$ the smallest constants $C$ such that the bounds on $\Pi, \Gamma$ and $\Sigma$ in the above analytic bounds hold. Furthermore, we define
$$\interleave Z\interleave_{\gamma;T}:= \|\Pi\|_{\gamma;T}+\|\Gamma\|_{\gamma;T}+\|\Sigma\|_{\gamma;T}.$$
If $\bar{Z}=(\bar{\Pi}, \bar{\Gamma},\bar{\Sigma})$ is another model we define
$$\interleave Z;\bar{Z}\interleave_{\gamma;T}:=\|\Pi-\bar{\Pi}\|_{\gamma;T}+\|\Gamma-\bar{\Gamma}\|_{\gamma;T}+\|\Sigma-\bar{\Sigma}\|_{\gamma;T},$$
This
gives a natural topology for the space of all models for a given regularity structure. In the following we consider the models are periodic in space, which allows us to require the bounds to hold globally.

Now we have the following definition for the spaces of distributions $\mathcal{C}^\alpha$, $\alpha<0$, which is an extension of the definition of H\"{o}lder space to include $\alpha<0$.
\vskip.10in

\th{Definition 2.3} For $\alpha<0$, $\mathcal{C}^\alpha$ consists of the closure of smooth compact functions under the norm
$$\|\eta\|_{\alpha}:=\sup_{x\in\mathbb{R}^d}\sup_{\varphi\in\mathcal{B}_r}\sup_{\lambda\leq 1}\lambda^{-\alpha}|\eta(\varphi_x^\lambda)|<\infty.$$

\vskip.10in

We also have the following definition of spaces of inhomogeneous modelled distributions, which are the  H\"{o}lder spaces on the regularity structure.

\vskip.10in
\th{Definition 2.4} Given a model $Z=(\Pi,\Gamma,\Sigma)$ for a regularity structure $\mathfrak{T}$   as above. Then for any $\gamma>0$ and $\eta\in\mathbb{R}$, the space $\mathcal{D}^{\gamma,\eta}$  consists of all functions $H:(0,T]\times\mathbb{R}^{d}\rightarrow \bigoplus_{\alpha<\gamma}T_\alpha$ such that
$$\interleave H\interleave_{\gamma,\eta;T}:=\|H\|_{\gamma,\eta;T}+\sup_{\tiny\aligned
s\neq t\in&(0,T]\\
|t-s|\leq& |t,s|_0^{\mathfrak{s}_0}
\endaligned }\sup_{x\in\mathbb{R}^d}
\sup_{l<\gamma}\frac{\|H_t(x)-\Sigma_x^{ts}H_s(x)\|_l}{|t-s|^{(\gamma-l)/\mathfrak{s}_0}|t,s|_0^{\eta-\gamma}}<\infty,$$
with
$$\| H\|_{\gamma,\eta;T}:=\sup_{t\in(0,T]}\sup_{x\in\mathbb{R}^d}\sup_{l<\gamma}|t|_0^{(l-\eta)\vee 0}\|H_t(x)\|_l+\sup_{t\in(0,T]}\sup_{\tiny\aligned
x\neq y\in&\mathbb{R}^d\\
|x-y|\leq&1
\endaligned }
\sup_{l<\gamma}\frac{\|H_t(x)-\Gamma_{xy}^{t}H_t(y)\|_l}{|x-y|^{\gamma-l}|t|_0^{\eta-\gamma}},$$
Here  we wrote $\|\tau\|_l$ for the norm of the component of $\tau$ in $T_l$ and $|t|_0:=|t|^{\frac{1}{\mathfrak{s}_0}}\wedge1$ and $|t,s|_0:=|t|_0\wedge|s|_0$.
\vskip.10in

For $H\in \mathcal{D}^{\gamma,\eta}$ and $\bar{H}\in \bar{\mathcal{D}}^{\gamma,\eta}$ (denoting by $\bar{\mathcal{D}}^{\gamma,\eta}$ the space built over another model $(\bar{\Pi},\bar{\Gamma},\bar{\Sigma})$), we also set
$$\aligned\| H;\bar{H}\|_{\gamma,\eta;T}:=&\sup_{t\in(0,T]}\sup_{x\in\mathbb{R}^d}\sup_{l<\gamma}|t|_0^{(l-\eta)\vee 0}\|H_t(x)-\bar{H}_t(x)\|_l\\&+\sup_{t\in(0,T]}\sup_{\tiny\aligned
x\neq y\in&\mathbb{R}^d\\
|x-y|\leq&1
\endaligned }
\sup_{l<\gamma}\frac{\|H_t(x)-\Gamma_{xy}^{t}H_t(y)-\bar{H}_t(x)+\bar{\Gamma}_{xy}^{t}\bar{H}_t(y)\|_l}{|x-y|^{\gamma-l}|t|_0^{\eta-\gamma}},\endaligned$$
$$\interleave H;\bar{H}\interleave_{\gamma,\eta;T}:=\|H;\bar{H}\|_{\gamma,\eta;T}+\sup_{\tiny\aligned
s\neq t\in&(0,T]\\
|t-s|\leq& |t,s|_0
\endaligned }\sup_{x\in\mathbb{R}^d}
\sup_{l<\gamma}\frac{\|H_t(x)-\Sigma_x^{ts}H_s(x)-\bar{H}_t(x)+\bar{\Sigma}^{ts}_x\bar{H}_s(x)\|_l}{|t-s|^{(\gamma-l)/\mathfrak{s}_0}|t,s|_0^{\eta-\gamma}},$$
which gives a natural
distance between elements $H\in \mathcal{D}^{\gamma,\eta}$ and $\bar{H}\in \bar{\mathcal{D}}^{\gamma,\eta}$.
\vskip.10in
Given a regularity structure, we say that a subspace $V\subset T$ is a sector of regularity $\alpha$ if it is invariant under the action of the structure group $G$ and it can be written as $V=\oplus_{\beta\in A}V_\beta$ with $V_\beta\subset T_\beta$, and $V_\beta=\{0\}$ for $\beta<\alpha$. We will use $\mathcal{D}^{\gamma,\eta}(V)$ to denote all functions in $\mathcal{D}^{\gamma,\eta}$ taking values in $V$.

The reconstruction theorem, which defines the so-called reconstruction operator, is one of the most fundamental result in
 the regularity structures theory.

\vskip.10in
\th{Theorem 2.5} (cf. [16, Theorem 2.11]) Given a model $Z=(\Pi,\Gamma,\Sigma)$ for a regularity structure $\mathfrak{T}$ with $\alpha:=\min A$ . Then for every $\eta\in\mathbb{R}, \gamma>0$ and $T>0$, there is a unique family of linear operators $\mathcal{R}_t:\mathcal{D}^{\gamma,\eta}\rightarrow \mathcal{C}^\alpha(\mathbb{R}^d)$, parametrised by $t\in(0,T]$, such that the bound
$$|\langle \mathcal{R}_tH_t-\Pi^t_xH_t(x),\varphi^\lambda_x\rangle|\lesssim \lambda^{\gamma}|t|_0^{\eta-\gamma}\|H\|_{\gamma,\eta;T}\|\Pi\|_{\gamma;T},$$
holds uniformly in $H\in\mathcal{D}^{\gamma,\eta}, t\in(0,T], x\in\mathbb{R}^d, \lambda\in (0,1]$ and $\varphi\in \mathcal{B}_r$ with $r>-\alpha+1$.

\vskip.10in

Suppose that $K$ is a $2$-regularising kernel in the sense of [14, Section 5] and we will write $K_t(x)=K(z)$, for $z=(t,x)$.
We say that a model $Z=(\Pi,\Gamma,\Sigma)$ realises $K$ for an abstract integration map $\mathcal{I}$ if, for every $\alpha\in A$, every $\tau\in T_\alpha$ and every $x\in\mathbb{R}^d$, one has \begin{equation}\Pi_x^t(\mathcal{I}\tau+\mathcal{J}_{t,x}\tau)(y)=\int_{\mathbb{R}}\langle \Pi^s_x\Sigma^{st}_x\tau, K_{t-s}(y-\cdot)\rangle ds.\end{equation}
Here $$\mathcal{J}_{t,x}\tau=\sum_{|k|_{\mathfrak{s}}<\alpha+\beta}\frac{X^k}{k!}\int_{\mathbb{R}}\langle \Pi^s_x\Sigma^{st}_x\tau, D^kK_{t-s}(x-\cdot)\rangle ds,$$
where $k\in\mathbb{N}^{d+1}$ and the derivative $D^k$ is in time-space. Moreover, we require that
\begin{equation}\Gamma^t_{xy}(\mathcal{I}+\mathcal{J}_{t,y})=(\mathcal{I}+\mathcal{J}_{t,x})\Gamma^t_{xy}, \quad \Sigma_x^{st}(\mathcal{I}+\mathcal{J}_{t,x})=(\mathcal{I}+\mathcal{J}_{s,x})\Sigma_x^{st},\end{equation}
for all $s,t\in\mathbb{R}$, and $x,y\in\mathbb{R}^d$.

Using $\mathcal{I}$ we can also introduce the  operator $\mathcal{K}_\gamma$ acting on modelled distribution similarly as in \cite{HM15} and \cite{ZZ15}.

\vskip.10in
In order to deal with the Leray Projection, we have to consider convolution with the singular kernel for space variable. As in \cite{Hai14} we introduce an abstract integration map $\mathcal{I}_0:T\rightarrow T$ to provide an "abstract" representation of the Leray Projection operating at the level of the regularity structure. In the regularity structure theory  $\mathcal{I}_0$ is a linear map from $T$ to $T$ such that $\mathcal{I}_0T_\alpha\subset T_{\alpha}$ and $\mathcal{I}_0\bar{T}=0$ and for every $\Gamma\in G, \tau\in T$ one has $\Gamma \mathcal{I}_0\tau-\mathcal{I}_0\Gamma\tau\in\bar{T}$.

We say that $P$ is a  $0$-regularising kernel on $\mathbb{R}^d$ if one can write $P=\sum_{n\geq0}P_n$, where each  $P_n:\mathbb{R}^{d}\rightarrow\mathbb{R}$ is smooth and compactly supported in a ball of radius $2^{-n}$ around the origin. Furthermore, we assume that for every multi-index $k$, one has a constant $C$ such that
$$\sup_x|D^kP_n(x)|\leq C2^{n(d+|k|)},$$
holds uniformly in $n$. Finally, we assume that $\int P_n(x)E(x)dz=0$ for every polynomial $E$ of degree at most $r$ for some sufficiently large value of $r$.
A model $Z=(\Pi,\Gamma,\Sigma)$ realises a $0$-regularising kernel $P$ on $\mathbb{R}^d$ for an abstract integration map $\mathcal{I}_0$ if, for every $\alpha\in A$, every $\tau\in T_\alpha$ and every $x\in\mathbb{R}^d$, one has \begin{equation}\Pi_x^t(\mathcal{I}_0\tau+\mathcal{J}^0_{t,x}\tau)(y)=\langle \Pi^t_x\tau, P(y-\cdot)\rangle,\end{equation}
where $$\mathcal{J}_{t,x}^0\tau=\sum_{|k|<\alpha}\frac{X^k}{k!}\langle \Pi^t_x\tau, D^kP(x-\cdot)\rangle ,$$
with $k\in\mathbb{N}^{d}$ and the derivative $D^k$ being the space derivative. Moreover, we require that
\begin{equation}\Gamma^t_{xy}(\mathcal{I}_0+\mathcal{J}_{t,y}^0)=(\mathcal{I}_0+\mathcal{J}_{t,x}^0)\Gamma^t_{xy}, \quad \Sigma_x^{st}(\mathcal{I}_0+\mathcal{J}_{t,x}^0)=(\mathcal{I}_0+\mathcal{J}_{s,x}^0)\Sigma_x^{st},\end{equation}
for all $s,t\in\mathbb{R}$, and $x,y\in\mathbb{R}^d$.

Now we introduce the following operator acting on modelled distribution $H\in\mathcal{D}^{\gamma,\eta}$ with $\gamma+\beta>0$:
$$(\mathcal{P}_{\gamma,t} H_t)(x):=\mathcal{I}_0H_t(x)+\mathcal{J}_{t,x}^0H_t(x)+(\mathcal{N}^0_{\gamma,t} H_t)(x).$$
Here
$$(\mathcal{N}^0_{\gamma,t} H_t)(x):=\sum_{|k|<\gamma}\frac{X^k}{k!}\langle \mathcal{R}_tH_t-\Pi^t_xH_t(x), D^kP(x-\cdot)\rangle ,$$
where $k\in\mathbb{N}^{d}$ and the derivative $D^k$ is in space. By [25, Theorem 2.7] we have $$\mathcal{R}_t(\mathcal{P}_{\gamma,t} H_t)(x)= \langle \mathcal{R}_tH_t,P(x-\cdot)\rangle .$$

\vskip.10in

In the following we extend [14, Thm 5.14] for the inhomogeneous model.
\vskip.10in

\th{Proposition 2.6} Let $\mathfrak{T}=(A,T,G)$ be a regularity structure containing the canonical regularity structure. Let $V\subset T$ be a sector of order $\bar{\gamma}$. Let $W\subset V$ be a subsector of $V$ and let $K$ be a $2$-regularising kernel and $P$ be a $0$-regularising kernel on $\mathbb{R}^d$. Let $(\Pi,\Gamma,\Sigma)$ be a model for $\mathfrak{T}$, and let $\mathcal{I}: W\rightarrow T$ be an abstract integration map of order $2$ such that $\Pi$ realises $K$ for $\mathcal{I}$ and let $\mathcal{I}_0: W\rightarrow T$ be an abstract integration map of order $0$ such that $\Pi$ realises $P$ for $\mathcal{I}_0$.

Then, there exists a regularity structure $\hat{\mathfrak{T}}$ containing $\mathfrak{T}$, a model $(\hat{\Pi},\hat{\Gamma},\hat{\Sigma})$ for $\hat{\mathfrak{T}}$ extending $(\Pi,\Gamma,\Sigma)$, and abstract integration maps $\hat{\mathcal{I}}$ of order $2$, $\hat{\mathcal{I}}_0$ of order $0$ acting on $\hat{V}=\iota V$ such that:
\begin{itemize}
  \item The model $\hat{\Pi}$ realises $K$ for $\hat{\mathcal{I}}$ and realises $P$ for $\hat{\mathcal{I}}_0$.
  \item The map $\hat{\mathcal{I}}$ and $\hat{\mathcal{I}}_0$ extend $\mathcal{I}$ and $\mathcal{I}_0$ in the sense that $\hat{\mathcal{I}}\iota a=\iota \mathcal{I} a, \hat{\mathcal{I}}_0\iota a=\iota \mathcal{I}_0 a$ for every $a\in W$.
\end{itemize}
Furthermore, the map $(\Pi,\Gamma,\Sigma)\mapsto (\hat{\Pi},\hat{\Gamma},\hat{\Sigma})$ is locally bounded and Lipschitz continuous.

\proof As in [14, Theorem 5.14] we can assume without loss of generality that the sector $V$ is given by a finite sum $V=V_{\alpha_1}\oplus V_{\alpha_2}\oplus...\oplus V_{\alpha_n}$, where the $\alpha_i$ are an increasing sequence of elements in $A$, and $W_{\alpha_k}=V_{\alpha_k}$ for $k<n$. We then denote by $\bar{W}$ the complement of $W_{\alpha_n}$ in $V_{\alpha_n}$ so that $V_{\alpha_n}=W_{\alpha_n}\oplus \bar{W}_{\alpha_n}$.  By similar arguments as in the proof of [14, Theorem 5.14] we can extend the regularity structure $\hat{T}=T\oplus \bar{W}$ with $\beta=2$ and define $\hat{\mathcal{I}}$ similarly.  For $a\in T, b\in \bar{W}$ we define $\hat{\Pi}^t_x$ to be given by
$$\hat{\Pi}^t_x(a,b)=\Pi_x^ta+\int_{\mathbb{R}}\langle \Pi^s_x\Sigma^{st}_xb(dy), K_{t-s}(\cdot-y)\rangle ds-\Pi^t_x\mathcal{J}_{t,x}b,$$
where ${J}_{t,x}$ is given by (2.1). By [14, Lemma 5.19] $\hat{\Pi}^t_x$ satisfies the required bounds when tested against smooth test functions that are localized near $x$.
We set
$$\hat{\Gamma}_{xy}^t=(\Gamma_{xy}^t,M_{xy}^t), \quad M_{xy}^tb=\mathcal{J}_{t,x}\Gamma_{xy}^tb-\Gamma_{xy}^t\mathcal{J}_{t,y}b.$$
$$\hat{\Sigma}_{x}^{ts}=(\Sigma_{x}^{ts},{M}_{x}^{ts}), \quad {M}_{x}^{ts}b=\mathcal{J}_{t,x}\Sigma_{x}^{ts}b-\Sigma_{x}^{ts}\mathcal{J}_{s,x}b.$$
The corresponding algebraic identity and the corresponding analytic bounds can be obtained by similar arguments as in the proof of  [14, Theorem 5.14].
Furthermore, similarly  we can extend the regularity structure  with $\beta=0$ and define $\hat{\mathcal{I}}_0$ similarly.  For $a\in \hat{T}, b\in \bar{W}$ we define $\hat{\hat{\Pi}}^t_x$ to be given by
$$\hat{\hat{\Pi}}^t_x(a,b)=\hat{\Pi}_x^ta+\langle \Pi^t_xb(dy), P(\cdot-y)\rangle -\Pi^t_x\mathcal{J}^0_{t,x}b,$$
where ${J}_{t,x}^0$ is given by (2.3). By [14, Lemma 5.19] $\hat{\hat{\Pi}}^t_x$ satisfies the required bounds when tested against smooth test functions that are localized near $x$.
We set
$$\hat{\hat{\Gamma}}_{xy}^t=(\hat{\Gamma}_{xy}^t,\hat{M}_{xy}^t), \quad \hat{M}_{xy}^tb=\mathcal{J}_{t,x}^0\Gamma_{xy}^tb-\Gamma_{xy}^t\mathcal{J}^0_{t,y}b.$$
$$\hat{\hat{\Sigma}}_{x}^{ts}=(\hat{\Sigma}_{x}^{ts},\hat{M}_{x}^{ts}), \quad \hat{M}_{x}^{ts}b=\mathcal{J}^0_{t,x}\Sigma_{x}^{ts}b-\Sigma_{x}^{ts}\mathcal{J}^0_{s,x}b.$$
 The corresponding algebraic identity and the analytic bounds for $\hat{\Gamma}_{xy}^t$, $\hat{M}_{xy}^tb$ and $\hat{\Sigma}_{x}^{ts}$ can be checked by similar arguments as in the proof of [14, Theorem 5.14]. To obtain the analytic bounds for $\hat{M}_{x}^{st}$, we introduce the following model for $T$:
$$(\tilde{\Pi}_{(t,x)}\tau)(s,y)=(\Pi^{s}_x\Sigma^{st}_x\tau)(y), \quad \tilde{\Gamma}_{(t,x),(s,y)}=\Gamma^{t}_{xy}\Sigma^{{ts}}_y=\Sigma^{{ts}}_x\Gamma^{s}_{xy},$$
which is a model in the original sense of [14, Def 2.17].  For $\tau\in T_l$,
  $k\in\mathbb{N}^3$, consider
$$((\mathcal{J}_{t,x}^{0}\Sigma^{ts}_x-\Sigma^{ts}_x\mathcal{J}_{s,x}^{0})\tau)_k,$$
with $\mathcal{J}^0$ defined in (2.3).
We decompose $\mathcal{J}^0$ as $\mathcal{J}^0=\sum_{n\geq0}\mathcal{J}^{0,(n)}_{t,x}$, where the nth term in each sum is obtained by replacing ${P}$ by ${P}_n$ in the expressions for $\mathcal{J}^0$.  Moreover, for $\tau\in T_l$
$$(\mathcal{J}^{0,(n)}_{t,x}\Sigma^{ts}_x\tau)_k=\frac{1}{k!}\sum_{|k|<\zeta< l}\langle \Pi^{t}_x\mathcal{Q}_\zeta\Sigma_x^{ts}\tau, D^k{P}_n(x-\cdot)\rangle,$$
$$(\Sigma^{ts}_x\mathcal{J}^{0,(n)}_{s,x}\tau)_k=\frac{1}{k!}\langle \Pi^{s}_x\tau, D^k{P}_n(x-\cdot)\rangle,$$
where $\mathcal{Q}_\zeta a$ denotes the component of $a$ in $T_\zeta$.
We first consider the case $2^{-n}\leq |t-s|^{\frac{1}{2}}$: by Definition 2.2 we have
$$|(\mathcal{J}^{0,(n)}_{t,x}\Sigma^{ts}_x\tau)_k|\lesssim \sum_{|k|<\zeta< l}2^{n|k|}2^{-n\zeta}|s-t|^{\frac{l-\zeta}{2}}\lesssim \sum_{\delta<0}|t-s|^{\frac{l-|k|+\delta}{2}}2^{\delta n},$$
and
$$|(\Sigma^{ts}_x\mathcal{J}^{0,(n)}_{s,x}\tau)_k|\lesssim 2^{n|k|}2^{-nl}\lesssim \sum_{\delta<0}|t-s|^{\frac{l-|k|+\delta}{2}}2^{\delta n}.$$
For the case that $ |t-s|^{\frac{1}{2}}\leq 2^{-n}$ we have
$$\aligned&(\mathcal{J}^{0,(n)}_{t,x}\Sigma^{ts}_x\tau)_k-(\Sigma^{ts}_x\mathcal{J}^{0,(n)}_{s,x}\tau)_k
\\=&-\frac{1}{k!}\sum_{\zeta\leq |k| }\langle \Pi^{t}_x\mathcal{Q}_\zeta\Sigma^{ts}_x\tau, D^k{P}_n(x-\cdot)\rangle+\frac{1}{k!}\langle\Pi^{t}_x\Sigma^{ts}_x\tau- \Pi^{s}_x\tau, D^k{P}_n(x-\cdot)\rangle\\=&:T_1^k+T_2^k,\endaligned$$
For $T_1^k$ we have
$$|T_1^k|\lesssim \sum_{\zeta< |k|}2^{n|k|}2^{-n\zeta}|s-t|^{\frac{l-\zeta}{2}}\lesssim\sum_{\delta>0}|t-s|^{\frac{l-|k|+\delta}{2}}2^{\delta n},$$
where the sum runs over a finite number of exponents.
In the following we consider $T_2^k$:
$$\aligned|T_2^k|=&\frac{1}{k!}|\langle \tilde{\Pi}_{(s,x)}\tau(t,\cdot)-\tilde{\Pi}_{(s,x)}\tau(s,\cdot), D^k{P}_n(x-\cdot)\rangle|\\=&\lim_{m\rightarrow\infty}\frac{1}{k!}|\sum_{(s_0,y)\in\Lambda_m^{\mathfrak{s}}}\langle \tilde{\Pi}_{(s,x)}\tau,\varphi^{m,\mathfrak{s}}_{(s_0,y)}\rangle\langle \varphi^{m,\mathfrak{s}}_{(s_0,y)}(t)-\varphi^{m,\mathfrak{s}}_{(s_0,y)}(s), D^kP_n(x-\cdot)\rangle|,\endaligned$$
where $\{\varphi^{m,\mathfrak{s}}_{(s_0,y)}\}\subset \mathcal{C}^r$ is the wavelet basis introduced in [14, Section 3], $\Lambda_m^\mathfrak{s}=\{\sum_{j=0}^d 2^{-m\mathfrak{s}_j}k_je_j:k_j\in\mathbb{Z}\},$
with $e_j$ denoting the jth element of the canonical basis of $\mathbb{R}^{d+1}$.
By the definition of the model we have
$$\aligned|\langle \tilde{\Pi}_{(s,x)}\tau,\varphi^{m,\mathfrak{s}}_{(s_0,y)}\rangle|=&|\langle \tilde{\Pi}_{(s_0,y)}\tilde{\Gamma}_{(s_0,y),(s,x)}\tau,\varphi^{m,\mathfrak{s}}_{(s_0,y)}\rangle|
\\\lesssim&\sum_{l_0<l}\|(s,x)-(s_0,y)\|_{\mathfrak{s}}^{l-l_0}2^{-\frac{m|\mathfrak{s}|}{2}-l_0m}.\endaligned$$
For $\langle \varphi^{m,\mathfrak{s}}_{(s_0,y)}(t)-\varphi^{m,\mathfrak{s}}_{(s_0,y)}(s), D^kP_n(x-\cdot)\rangle$ we choose $m$ large enough such that $2^{-m}\leq |t-s|^{\frac{1}{\mathfrak{s}_0}}\leq 2^{-n}$. In this case by a similar calculation as in the proof of [14, Theorem 3.10] we know that
$$|\langle \varphi^{m,\mathfrak{s}}_{(s_0,y)}(t)-\varphi^{m,\mathfrak{s}}_{(s_0,y)}(s), D^kP_n(x-\cdot)\rangle|\lesssim 2^{n|k|}2^{-\frac{3m}{2}-rm}2^{n(3+r)}2^{\frac{\mathfrak{s}_0}{2}m}.$$
Furthermore, $|\langle \varphi^{m,\mathfrak{s}}_{(s_0,y)}(t)-\varphi^{m,\mathfrak{s}}_{(s_0,y)}(s), D^kP_n(x-\cdot)\rangle|=0$
unless $|x-y|\lesssim 2^{-n}$ and $|s-s_0|^{\frac{1}{\mathfrak{s}_0}}\wedge|t-s_0|^{\frac{1}{\mathfrak{s}_0}}\lesssim 2^{-m}$. Thus we have
$$\aligned|T_2^k|\lesssim&\lim_{m\rightarrow\infty}\sum_{l_0<l}2^{3m}2^{-3n}
2^{-\frac{m|\mathfrak{s}|}{2}-l_0 m}2^{-n(l-l_0)}2^{n|k|}2^{-\frac{3m}{2}-rm}2^{n(3+r)}2^m
\\\lesssim&\sum_{l_0<l}\lim_{m\rightarrow\infty}2^{-l m}2^{n|k|}2^{(l-l_0-r)(m-n)}\\\lesssim&\sum_{\delta>0}|t-s|^{\frac{l-|k|+\delta}{2}}2^{\delta n},\endaligned$$
where the sum runs over a finite number of exponents and in the first inequality we used  the factor $2^{3m}2^{-3n}$ counts the number of non-zero terms appearing in the sum over $(s_0,y)$ and in the last inequality we choose $m$ large enough such that $2^{-m}<|t-s|^{\frac{1}{2}}$ and $r$ large enough such that $r>l-l_0$.
Taking the sum over $n$ we obtain the desired bounds for $\|\hat{M}^{ts}_x\tau\|_k$ and the result follows.$\hfill\Box$

\subsection{N-S equation on $\mathbb{T}^3$}

In this subsection we recall the regularity structure theory for the 3D Navier-Stokes equations on $\mathbb{T}^3$ driven by space-time white noise in \cite{ZZ15}.
In this case we have the scaling $\mathfrak{s}=(2,1,1,1)$, so that the scaling dimension of space-time is $5$.
Since  the heat kernel $G$ is smooth on $\mathbb{R}^{4}\backslash\{0\}$ and has the scaling property $G(\frac{t}{\delta^2},\frac{x}{\delta})=\delta^3G(t,x)$ for $\delta>0$, by [14, Lemma 5.5] it can be decomposed into $K+R$ where $K$ is a $2$-regularising kernel and $R\in \mathcal{C}^\infty$

We know that the kernel $P^{ij}, i,j=1,2,3,$ for the Leray projection  is smooth on $\mathbb{R}^{3}\backslash\{0\}$ and has the scaling property $P^{ij}(\frac{x}{\delta})=\delta^3P^{ij}(x)$ for $\delta>0$, by [14, Lemma 5.5] it can be decomposed into $\bar{P}^{ij}+R_0^{ij}, i,j=1,2,3,$ where $\bar{P}^{ij}$ is a $0$-regularising kernel on $\mathbb{R}^3$ and $R_0^{ij}\in \mathcal{C}^\infty$.
Define $$K^{ij}:=K*\bar{P}^{ij}.$$

Consider the regularity structure generated by the stochastic N-S equation with $\beta=2$. In the regularity structure we use symbol the $\Xi_i$ to replace the driving noise $\xi^i$. We introduce the integration map $\mathcal{I}$  associated with $K$ and the integration map $\mathcal{I}_0^{ij}$ associated with $\bar{P}^{ij}$, which helps us to define $\mathcal{K}_\gamma$ and $\bar{\mathcal{P}}_{\gamma}^{ij}$. We also need
the integration maps $\mathcal{I}^{ii_1}_{0,k}, i,i_1=1,2,3, \mathcal{I}_k$ for a multiindex $k$, which  represents integration against $D^kP^{ii_1}, i,i_1=1,2,3, D^k K$ respectively.

To apply the regularity structure theory we write the equation as follows: for $i=1,2,3$
\begin{equation}\aligned\partial_t v_1^i=&\nu \Delta v_1^{i}+\sum_{i_1=1}^3P^{ii_1}\xi^{i_1},\quad \textrm{div} v_1=0,\\
\partial_t v^i=&\nu \Delta v^{i}-\sum_{i_1,j=1}^3P^{ii_1}\frac{1}{2}D_j [(v^{i_1}+v_1^{i_1}) (v^j+v_1^j)],\quad \textrm{div} v=0.\endaligned\end{equation}
Then $v_1+v$ is the solution to the 3D Navier-Stokes equations driven by space-time white noise.
Now we consider the second equation in (2.5). Define for $i,j,i_1=1,2,3$,
$$\mathcal{I}^{ij}:=\mathcal{I}_0^{ij}\mathcal{I},\quad \mathcal{I}^{ij}_{i_1}:=\mathcal{I}_0^{ij}\mathcal{I}_{i_1}$$
$$\mathfrak{M}_F^{ij}=\{1,\mathcal{I}^{ii_1}(\Xi_{i_1}),\mathcal{I}^{jj_1}(\Xi_{j_1}), \mathcal{I}^{ii_1}(\Xi_{i_1})\mathcal{I}^{jj_1}(\Xi_{j_1}), U_i, U_j, U_iU_j, \mathcal{I}^{ii_1}(\Xi_{i_1})U_j, U_i\mathcal{I}^{jj_1}(\Xi_{j_1}), i_1,j_1=1,2,3\},$$
where the product is commutative and associative.
Then we build subsets $\{\mathcal{P}_n^i\}_{n\geq0}$, $\{\mathcal{U}_n\}_{n\geq0}$ and $\{\mathcal{W}_n\}_{n\geq0}$ by the following algorithm: For $i,j=1,2,3$, set $\mathcal{W}_0^{ij}=\mathcal{P}_0^i=\mathcal{U}_0=\varnothing$ and
$$\mathcal{W}_n^{ij}=\mathcal{W}_{n-1}^{ij}\cup\bigcup_{\mathcal{Q}\in \mathfrak{M}_F^{ij}}\mathcal{Q}(\mathcal{P}_{n-1}^i,\mathcal{P}_{n-1}^j),$$
$$\mathcal{P}_{n}^i=\{X^k\}\cup\{\mathcal{I}^{ii_1}_{i_2}(\tau):\tau\in\mathcal{W}_{n-1}^{i_1i_2}, i_1,i_2=1,2,3\},$$
$$\mathcal{U}_{n}=\{  \mathcal{I}_{i_2}(\tau):\tau\in\mathcal{W}_{n-1}^{i_1i_2},  i_1,i_2=1,2,3\},$$
and
$$\mathcal{F}_F:=\bigcup_{n\geq0}(\bigcup_{i,j=1}^3\mathcal{W}_n^{ij}\cup\mathcal{U}_{n} ).$$
Then $\mathcal{F}_F$ contains the elements required to describe both the solution and the terms  in the equation (2.5).
We denote by $\mathcal{H}_F$ the set of finite linear combinations of elements in $\mathcal{F}_F$. For each $\tau\in \mathcal{F}$ a degree $|\tau|_{\mathfrak{s}}$  is obtained by setting
$|\mathbf{1}|_\mathfrak{s}=0$,
$$|\tau\bar{\tau}|_\mathfrak{s}=|\tau|_\mathfrak{s}+|\bar{\tau}|_\mathfrak{s},$$
for any two formal expressions $\tau$ and $\bar{\tau}$ in $\mathcal{F}$ such that
$$|\Xi_i|_\mathfrak{s}=\alpha,\quad |X_i|_\mathfrak{s}=\mathfrak{s}_i,\quad|\mathcal{I}_{k}(\tau)|_\mathfrak{s}=|\tau|_\mathfrak{s}+2-|k|_\mathfrak{s}, \quad |\mathcal{I}^{ii_1}_{0,k}(\tau)|_\mathfrak{s}=|\tau|_\mathfrak{s}-|k|_\mathfrak{s},$$
for $\alpha\in (-\frac{13}{5},-\frac{5}{2})$. Let $T=\mathcal{H}_F$ with $T_\gamma=\langle \{\tau\in\mathcal{F}_F:|\tau|_\mathfrak{s}=\gamma\}\rangle$, $A=\{|\tau|_\mathfrak{s}:\tau\in\mathcal{F}_F\}$ and let the structure group $G_F$ be as in \cite{ZZ15}. Then by [25, Theorem 2.8] we know that $\mathfrak{T}_F=(A,\mathcal{H}_F, G_F)$ defines a regularity structure $\mathfrak{T}$. We also recall the following definition from [25, Definition 2.10].

\vskip.10in
\th{Definition 2.7} A model $(\Pi,\Gamma,\Sigma)$ for $\mathfrak{T}$ is admissible if it satisfies  \begin{equation}(\Sigma_x^{st}\tau\bar{\tau})=(\Sigma_x^{st}\tau)(\Sigma_x^{st}\bar{\tau}),\quad (\Gamma_{xy}^{t}\tau\bar{\tau})=(\Gamma_{xy}^{t}\tau)(\Gamma_{xy}^{t}\bar{\tau}),\end{equation} and furthermore realizes $K, \bar{P}^{ij}, i,j=1,2,3,$ for $\mathcal{I}, \mathcal{I}_0^{ij}$ respectively. We denote by $\mathcal{M}_F$ the set of admissible models.
\vskip.10in

Using the same tree notation from \cite{ZZ15} we can denote $\mathcal{F}_0:=\mathcal{T}_{<0}$ as follows:

For $\Xi$ we simply draw a dot. The integration map $\mathcal{I}^{ij}$ is then represented by a downfacing line while the integration map  $\mathcal{I}_0\mathcal{I}_j$ is then represented by a downfacing dotted line. The integration map $\mathcal{I}_j$ is represented by $\includegraphics[height=0.3cm]{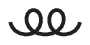}$. The multiplication of symbols is obtained by joining them at the root.
$$\aligned\mathcal{F}_0=\{\mathbf{1}, &\includegraphics[height=0.5cm]{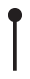}, \includegraphics[height=0.5cm]{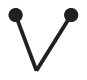},\includegraphics[height=0.7cm]{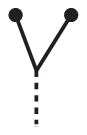}, \includegraphics[height=0.7cm]{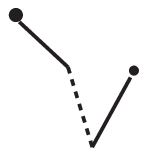},\includegraphics[height=0.7cm]{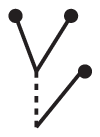},
\includegraphics[height=0.7cm]{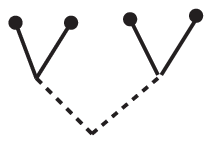},\includegraphics[height=0.7cm]{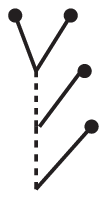}, \includegraphics[height=0.7cm]{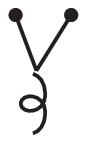}\},\endaligned$$
Here the solution to the  stochastic N-S equation is vector valued and there are a lot of superscripts and subscripts for the elements in $\mathcal{F}_0$, which  will not be noticeable in the tree notation. To see more details for $\mathcal{F}_0$ we refer to \cite{ZZ15}.

Similarly as in \cite{Hai14}, the renormalization group $\mathfrak{R}$ for the regularity structure $\mathfrak{T}_F$ is also introduced in \cite{ZZ15}.
For $g \in \mathfrak{R}:=\mathbb{R}^n, n=3^4+3\cdot3^{10}$, $g=(C^1_{ii_1jj_1}, C^2_{ii_1i_2jj_1j_2kk_1ll_1}, C^3_{ii_1i_2i_3kk_1ll_1jj_1}, C^4_{ii_1i_2kk_1ll_1l_2jj_1}),$  $i,j,k,l, i_1,i_2,i_3,j_1,$ $k_1,l_1,l_2=1,2,3$, we define a linear map $M_g$ on $\mathcal{F}_0$ by
$$\aligned M_g\includegraphics[height=0.5cm]{02.eps}=&\includegraphics[height=0.5cm]{02.eps}-C^1_{ii_1jj_1}\mathbf{1},
\\M_g\includegraphics[height=0.7cm]{07.eps}=&\includegraphics[height=0.7cm]{07.eps}-C^2_{ii_1i_2jj_1j_2kk_1ll_1}\mathbf{1},
\\M_g\includegraphics[height=0.7cm]{08.eps}=&\includegraphics[height=0.7cm]{08.eps}-C^p_{ii_1i_2i_3kk_1ll_1jj_1}\mathbf{1},\quad p=3,4\endaligned$$
as well as $M_g(\tau)=\tau$ for the remaining basis vectors in $\mathcal{F}_0$. Here the choice that $p=3$ or $4$ depends on the explicit formula of $\includegraphics[height=0.7cm]{08.eps}$. For more details, we refer to \cite{ZZ15}.

For given $\xi\in \mathcal{C}^\infty(\mathbb{R}\times \mathbb{R}^3;\mathbb{R}^3)$ by [25, Section 2] we can construct a canonical model $\mathcal{L}\xi\in\mathcal{M}_F$. For $g\in\mathfrak{R}$, we can define $M_g$ on $\mathcal{F}_0$ as above. Since $\tau$ satisfies $M_g\tau=\tau-C\mathbf{1}$ for any $\tau\in\mathcal{F}_0$, we can easily lift this action $M_g$ on the space $\mathcal{M}_F$ of admissible models via the construction of [14, Section 8] and [25, Section 2]. We also use the following definition from \cite{HM16}. We choose  $\eta\in (-1,\alpha+2]$ and $|\alpha+2|<\gamma<\eta-\alpha$ for $\alpha\in (-\frac{13}{5},-\frac{5}{2})$. Here $\alpha$ is the degree of $\Xi_i$ and the initial value belongs to $\mathcal{C}^\eta$.
\vskip.10in

\th{Definition 2.8} An admissible model $\mathbf{\Pi}=(\Pi,\Gamma,\Sigma)\in\mathcal{M}_F$ is \emph{nice} if there exist $\xi_n\in\mathcal{C}^\infty(\mathbb{R}\times \mathbb{R}^d)$ and $g_n\in\mathfrak{R}$ such that $\mathbf{\Pi}=\lim_{n\rightarrow} M_{g_n}\mathcal{L}(\xi_n)$ and furthermore the distribution $K^{ii_1}*\xi^{i_1}:=\mathcal{R}_\cdot\mathcal{I}^{ii_1}(\Xi_{i_1})\in C(\mathbb{R},\mathcal{C}^{\eta})$, $D_jK^{i_0i}*(K^{ii_1}*\xi^{i_1}\diamond K^{jj_1}*\xi^{j_1})\in C(\mathbb{R},\mathcal{C}^{2\alpha+5}) $, $i,i_0,i_1,j,j_1=1,2,3$. Here $\mathcal{R}_t$ is the family of reconstruction operators in Theorem 2.5.
\vskip.10in

Write $\mathcal{M}$ for the closure of all smooth and nice models in the space of nice admissible models for the regularity structure $\mathfrak{T}_F$.
Now we  consider the following equations on the regularity structure:
 \begin{equation}\aligned u^i=&-\frac{1}{2}\sum_{i_1,j=1}^3\mathcal{T}^{ii_1}_j\mathbf{R}^+ (u^{i_1}\star u^j)+\sum_{i_1=1}^3\bar{\mathcal{T}}^{ii_1}_j\mathbf{R}^+\Xi_{i_1}+\sum_{i_1=1}^3\mathcal{G}^{ii_1}u_0^{i_1}
 \\J^i=&-\frac{1}{2}\sum_{i_1,j=1}^3\mathcal{T}^{ii_1}_j\mathbf{R}^+ (u^{i_1}\star J^j+J^{i_1}\star u^j)+\sum_{i_1=1}^3\mathcal{G}^{ii_1}J_0^{i_1}
,\endaligned\end{equation}
with the convolution operators $\mathcal{T}^{ii_1}_j$   and   $\bar{\mathcal{T}}^{ii_1}_j$ satisfying $$\mathcal{R}_t(\mathcal{T}^{ii_1}_jH)_t(x)=\int_0^t \langle (P^{ii_1}*D_jG)_{t-s}(x-\cdot),\mathcal{R}_sH_s\rangle ds,$$
  $$\mathcal{R}_t(\bar{\mathcal{T}}^{ii_1}_jH)_t(x)=\int_0^t \langle (P^{ii_1}*G)_{t-s}(x-\cdot),\mathcal{R}_sH_s\rangle ds,$$
  for $H\in \mathcal{D}^{\gamma,\eta}$,  and $\mathcal{G}^{ii_1}$ denotes the solution to the linearized problem, viewed via its truncated Taylor expansion as an element in $\mathcal{D}^{\gamma,\eta}$. Here $\mathbf{R}^+:\mathbb{R}\times \mathbb{R}^3\rightarrow\mathbb{R}$ is given by $\mathbf{R}^+(t,x)=1$ for $t>0$ and $\mathbf{R}^+(t,x)=0$ otherwise and for more details on the convolution operators $\mathcal{T}^{ii_1}_j$   and   $\bar{\mathcal{T}}^{ii_1}_j$, we refer to \cite{ZZ15}. By similar arguments as in the proof of [25, Proposition 2.13], we have the following results:

\vskip.10in

\th{Proposition 2.9} Let $\mathfrak{T}_F$ be the regularity structure constructed above associated to the stochastic N-S equation driven by space-time white noise. Suppose that $\eta\in (-1,\alpha+2]$ and $|\alpha+2|<\gamma<\eta-\alpha$ for $\alpha\in (-\frac{13}{5},-\frac{5}{2})$. Let $u_0^i, J_0^i\in \mathcal{C}^\eta(\mathbb{R}^3), i=1,2,3$, periodic and let $\mathbf{\Pi}=(\Pi,\Gamma,\Sigma)\in\mathcal{M}$ be a nice model for $\mathfrak{T}_F$. Then there exists a maximal solution  to  equations (2.7).

\proof This result for the first equation in (2.7) has been obtained in [25, Proposition 2.13]. It is sufficient to consider the second one in (2.7). We have that $J$ takes values in a sector $V\subset \bar{T}\oplus T_{\geq \zeta}$ with $\zeta=\alpha+3$ and $u^iJ^j, i,j=1,2,3,$ takes value in a sector $V\subset \bar{T}\oplus T_{\geq \bar{\zeta}}$ with  $\bar{\zeta}=2\alpha+5$ satisfying $\zeta<\bar{\zeta}+1$. For $\eta$ and $\gamma$ we have $\bar{\eta}=2\eta$ and $\gamma>\bar{\gamma}=\gamma+\alpha+2>0$, $\bar{\gamma}<\eta+2$ and  $\bar{\gamma}+1>\gamma$. Then by similar arguments as in the proof of [25, Proposition 2.13] the results follow.$\hfill\Box$

\section{Strong Feller Property}

We consider a Gaussian probability space $(\Omega,\mathcal{F},\mathbf{P})$, where $\Omega$ is a separable Banach space and $L^2(\mathbb{R};L^2)$ is the Cameron-Martin space for the Gaussian measure $\mathbf{P}$. The canonical random variable $\omega$ induces the usual two-sided continuous filtration $\{\mathcal{F}_{s,t},s<t\}$. $\mathcal{F}_t:=\mathcal{F}_{-\infty,t}$. Let $\xi= (\xi^1, \xi^2, \xi^3)$, with $\xi^i, i = 1, 2, 3$ being independent white noises
on $\mathbb{R}\times \mathbb{T}^3$  on $(\Omega,\mathcal{F},\mathbf{P})$, which we extend periodically to $\mathbb{R}^4$. By [25, Theorem 2.17] we know that there exists a random variable $\mathbf{Z}:\Omega \rightarrow\mathcal{M}$ associated with space-time white noise $\xi$.

Choose $U=(\mathcal{C}^\eta)^3$ with the usual product norm still denoted by $\|\cdot\|_\eta$ and define $\bar{U}=U\cup \{\Delta\}$. $\bar{U}$ is a separable metric space by setting $d(\Delta,u)^2=1+\|u\|_\eta^2$ for all $u\in U$. For $u_0\in U, \mathbf{\Pi}\in \mathcal{M}$  we use $\bar{u}_{s,t}(u_0,\mathbf{\Pi})$ to denote the maximal solution to the first equation in (2.7) obtained in Proposition 2.9 with initial condition $\bar{u}_{s,s}(u_0,\mathbf{\Pi})=u_0\in U$. Define
$${u}_{s,t}(u_0,\mathbf{\Pi})=\left\{\begin{array}{ll}\mathcal{R}_t\bar{u}_{s,t}(u_0,\mathbf{\Pi}),&\ \ \ \ \textrm{ for } t\in [s,\zeta)
\\\Delta&\ \ \ \ \textrm{ for } t\in [\zeta,\infty),\end{array}\right.$$
and ${u}_{s,t}(\Delta,\mathbf{\Pi})=\Delta$, where $\zeta$ is the explosion time for $s, u_0$ and $\mathbf{\Pi}$. We also use the shorthands $u_t=u_{0,t}$ and $u=u_{0,1}$. It follows from  the locality of the
reconstruction map and the locality of the construction of the model that ${u}_{s,t}(u_0,\mathbf{\Pi})$
depends on the underlying white noise only on the time interval $[s,t]$. Moreover,
as a consequence of [14, Prop. 7.11], one has for $u_0\in U, \mathbf{\Pi}\in\mathcal{M}$
$${u}_{0,s+t}(u_0,\mathbf{\Pi})={u}_{s,s+t}(u_{0,s}(u_0,\mathbf{\Pi}),\mathbf{\Pi}).$$
By this we can conclude that $u_{\cdot}(u_0,\mathbf{Z})$ is a Markov process.
 For any bounded measurable map $\Psi:\bar{U}\rightarrow\mathbf{R}$ we set for $t\geq0$
 $$P_t\Psi(u_0):=\mathbf{E}\Psi(u_t(u_0,\mathbf{Z})),$$
 which forms a Markov semigroup.

\vskip.10in

\th{Theorem 3.1} The Markov semigroup $(P_t)_{t\geq0}$ satisfies the strong Feller property, i.e. $P_t\Psi$ is continuous for all $\Psi$ that are only bounded and measurable.

\proof It is sufficient to check that Assumptions 1-5 in \cite{HM16} hold. Here for  the reader's convenience, we summarize them in appendix. We first check Assumption 1: By [25, Propositions 2.13, 2.15] we know that the preimage of $U$ under the map $(s,t,u_0,\mathbf{\Pi})\mapsto u_{s,t}(u_0,\mathbf{\Pi})$ is open and the map is jointly continuous on the preimage of $U$.
 The continuity in time for $u$ follows from [25, Theorem 1.1]. Now we consider Frech\'{e}t differentiability with respect to the initial condition.
By   the implicit function theorem in [7, Theorem 19.28] and similar arguments as in the proof of [3, Proposition 4.7], it follows that the solution $u_t$ is differentiable in the initial condition and its derivative in the direction $J_0\in (\mathcal{C}^\eta)^3$ is given by $J$. Here $J_t=\mathcal{R}_t\bar{J}_t$ with $\bar{J}$ being the maximal solution to the second equation in (2.7) obtained in Proposition 2.9.

Now we check Assumption 2 in \cite{HM16}.
We set
$$r_t(u_0,\mathbf{\Pi})=\left\{\begin{array}{ll}+\infty,&\ \ \ \ \textrm{ if } u_t(u_0,\mathbf{\Pi})=\Delta
\\\interleave \bar{u}\interleave_{\gamma,\eta;t}&\ \ \ \ \textrm{otherwise}.\end{array}\right.$$
The $\mathcal{F}_t$-measurability of $r_t$ is an immediate consequence of the adaptedness of $u$. It is obvious that  $r_0<\infty$ and that the map $t\mapsto r_t(u_0,\mathbf{\Pi})$ is increasing and is continuous in $t$, except at the explosion time when it has to diverge to $+\infty$. By [25, Proposition 2.13] we know that $r$ is locally Lipschitz continuous on $\mathcal{N}_t$ as a function of both the initial condition and the underlying model.

Now we check Assumption 3 in \cite{HM16}. With $M_g$ and $\mathfrak{L}$ given in Section 2, we first check Assumption 10 in \cite{HM16} as in [17, Section 5.1]. For $p>1$ sufficiently large, we choose $E:=L^p([0,1],\mathcal{C}^3)$ as a suitable space of shifts with $\mathcal{C}$ the space of periodic continuous functions. We prove that $E$ generates a continuous action on our space $\mathcal{M}$ of nice models. We introduce an auxiliary regularity structure $(\hat{{T}},\hat{{G}})$ as follows.
We introduce $\hat{\Xi}_i, i=1,2,3$, for the shifts. Define for $i,j=1,2,3$,
$$\aligned\hat{\mathfrak{M}}_F^{ij}=&\{1,\mathcal{I}^{ii_1}(\bar{\Xi}_{i_1}),\mathcal{I}^{jj_1}(\bar{\Xi}_{j_1}), \mathcal{I}^{ii_1}(\bar{\Xi}_{i_1})\mathcal{I}^{jj_1}(\bar{\Xi}_{j_1}), U_i, U_j, U_iU_j, \mathcal{I}^{ii_1}(\bar{\Xi}_{i_1})U_j, U_i\mathcal{I}^{jj_1}(\bar{\Xi}_{j_1}),\\& \bar{\Xi}_{i_1}=\Xi_{i_1} \textrm{ or } \hat{\Xi}_{i_1}, \bar{\Xi}_{j_1}=\Xi_{j_1} \textrm{ or } \hat{\Xi}_{j_1} i_1,j_1=1,2,3\}.\endaligned$$
Then we build
$\hat{\mathcal{F}_F}$ by a similar argument as that for $\mathcal{F}_F$, but with ${\mathfrak{M}}_F^{ij}$ replaced by $\hat{\mathfrak{M}}_F^{ij}$.
Graphically, if we denote $\hat{\Xi}$ be a circle, the symbols appearing in $\hat{\mathcal{F}_F}$ are the same as those appearing in  $\mathcal{F}_F$
but with any occurrence of a bullet possibly replaced by a circle.

We set $|\hat{\Xi}|_{\mathfrak{s}}=-\kappa$ for $\kappa$ sufficiently small, the degrees of the remaining basis vectors being obtained by using the same rules as above.
The structure group $\hat{{G}}$ is also defined similarly as in \cite{ZZ15} and by imposing that $\Gamma\hat{\Xi}_j=\hat{\Xi}_j$ for $j=1,2,3$.
 Let $\hat{{T}}=\langle\hat{\mathcal{F}_F}\rangle$ with $\hat{T}_\gamma=\langle \{\tau\in\hat{\mathcal{F}}_F:|\tau|_\mathfrak{s}=\gamma\}\rangle$, $\hat{A}=\{|\tau|_\mathfrak{s}:\tau\in\hat{\mathcal{F}}_F\}$ and let $\hat{{G}}$ be as above. Then $\hat{\mathfrak{T}}_F=(\hat{A},\hat{{T}}, \hat{{G}})$ defines a regularity structure $\hat{\mathfrak{T}}$.

Define $\hat{\mathcal{F}}_0=\mathcal{F}_0\cup \{\mathcal{I}^{ii_1}(\hat{\Xi}_{i_1}), i,i_1=1,2,3\}$.
For   $g\in \mathfrak{R}$, we define a linear map $\hat{M}_g$ on $\langle\hat{\mathcal{F}}_0\rangle$ as in Section 2.2 and for $\tau\in \hat{\mathcal{F}_0}\setminus\mathcal{F}_0$, $\hat{M}_g\tau=\tau$.
As a consequence of the calculation in \cite{ZZ15}, $\hat{M}_g$ belongs to the renormalisation group  defined in [14,
Definition 8.43]. Then similarly as in [14, Theorem 8.46] we can define the related renormlised model $(\hat{\Pi}^{\hat{M}_g}, \hat{\Gamma}^{\hat{M}_g},\hat{\Sigma}^{\hat{M}_g})$ and it is
an admissible model for $\hat{\mathfrak{T}}_F$ on $\langle\hat{\mathcal{F}}_0\rangle$.
Furthermore, by Proposition 2.6 and [14, Prop. 3.31] it extends uniquely to an admissible model
for all of $\hat{\mathfrak{T}}_F$.  Furthermore,  $\mathcal{I}$ is an  abstract integration map of order $2$ on $\langle\hat{\mathcal{F}}_0\rangle$ and for every $i,i_1=1,2,3$, $\mathcal{I}_0^{ii_1}$ is an abstract integration map of order $0$ on $\langle\hat{\mathcal{F}}_0\rangle$. We also have a space of "nice models" $\hat{\mathcal{M}}$ for this large regularity structure $\hat{\mathfrak{T}}_F$.

Since $\kappa$ is sufficiently small, all of the elements of $\hat{\mathcal{F}}_F\backslash\mathcal{F}_F$ are of strictly positive degree. By repeatedly applying
Proposition 2.6 and writing $\hat{\mathcal{M}}_F$ as the space of admissible models for $(\hat{{T}},\hat{{G}})$, there exists a unique locally Lipschitz continuous map $\mathcal{Y}:E\times {\mathcal{M}}_F\rightarrow\hat{\mathcal{M}}_F$ such that
for every $h\in E$ and $\mathbf{\Pi}\in \mathcal{M}_F$, the model $\hat{\mathbf{\Pi}}=\mathcal{Y}(h,\mathbf{\Pi})\in \hat{\mathcal{M}}_F$ agrees with $\mathbf{\Pi}$ on ${\mathfrak{T}}_F$ and for $(\hat{\Pi}, \hat{\Gamma},\hat{\Sigma}):=\mathcal{Y}(h,\mathbf{\Pi})$ we have $$\hat{\Pi}_x^t\mathcal{I}^{ii_1}(\hat{\Xi}_{i_1})(y)=K^{ii_1}*h^{i_1}(t,y)-K^{ii_1}*h^{i_1}(t,x)-\sum_j(y_j-x_j)D_jK^{ii_1}*h^{i_1}(t,x).$$
Since $\hat{M}_g\mathcal{Y}(h,\mathbf{\Pi})$ and $\mathcal{Y}(h,M_g\mathbf{\Pi})$ agree on $\hat{\mathcal{F}}_0$ and any admissible model is uniquely determined by this, we have for every $g\in\mathfrak{R}$, every $\mathbf{\Pi}\in\mathcal{M}$ and every $h\in E$, $\hat{M}_g\mathcal{Y}(h,\mathbf{\Pi})=\mathcal{Y}(h,M_g\mathbf{\Pi})$.
Now as in \cite{HM16} we introduce a map $\mathcal{Z}:\hat{\mathcal{M}}_F\rightarrow\mathcal{M}_F$, also commuting with the action of $\mathfrak{R}$, so that we can define $\tau=\mathcal{Z}\circ\mathcal{Y}$. For this, given a model $\hat{\mathbf{\Pi}}\in \hat{\mathcal{M}}_F$ we introduce the following notion of a $\hat{\mathbf{\Pi}}$-polynomial as in \cite{HM16}:
\vskip.10in

\th{Definition 3.2} A $\hat{\mathbf{\Pi}}$-polynomial $f$ is a map $f:\mathbb{R}^{d+1}\rightarrow\hat{\mathcal{M}}_F$ such that $f(t,x)=\Gamma^t_{xy}f(t,y)$ and $f(t,x)=\Sigma^{ts}_xf(s,x)$ for every $s,t,x, y$.

\vskip.10in
Given an admissible model $\hat{\mathbf{\Pi}}=(\Pi,\Gamma,\Sigma)\in\hat{\mathcal{M}}_F$ we define a collection of $\hat{\mathbf{\Pi}}$-polynomial $\{f^\tau_{t,x}\}$ for $\tau\in \mathcal{F}_F$ and $(t,x)\in\mathbb{R}^{d+1}$ recursively as follows. Define deg$\tau:=|\tau|_{\mathfrak{s}}$ for $\tau\in \mathcal{F}_F$, but then setting deg$\hat{\Xi}:=|\Xi|_{\mathfrak{s}}$ and defining it on the rest of $\hat{\mathcal{F}}_F$ by using the same rules as for $|\cdot|_{\mathfrak{s}}$. This allows us to define operators $\bar{\mathcal{J}}_{t,x}:\hat{{T}}\rightarrow\bar{{T}}, \bar{\mathcal{J}}_{t,x}^{0,ij}:\hat{{T}}\rightarrow\bar{{T}}$ by setting
$$\bar{\mathcal{J}}_{t,x}\tau=\mathcal{Q}_{<\mathrm{deg}\tau+2}\mathcal{J}_{t,x}\tau, \quad \bar{\mathcal{J}}^{0,ij}_{t,x}\tau=\mathcal{Q}_{<\mathrm{deg}\tau}\mathcal{J}^{0,ij}_{t,x}\tau,$$
where $\mathcal{J}^{0,ij}_{t,x}$ is defined in (2.3) associated with $\bar{P}^{ij}$.
With these definitions at hand we set for $z=(t,x), \bar{z}=(\bar{t},\bar{x})$
$$f^{X^k}_z(\bar{z})=\Gamma^{\bar{t}}_{\bar{x}x}\Sigma_x^{\bar{t}t}X^k,\quad f_z^{\Xi_{j}}(\bar{z})=\Xi_{j}+\hat{\Xi}_{j}.$$
Then we set recursively
$$f_z^{\tau\bar{\tau}}(\bar{z})=f_z^\tau(\bar{z})f_z^{\bar{\tau}}(\bar{z}),\quad f_z^{\mathcal{I}(\tau)}(\bar{z})=(\mathcal{I}+\mathcal{J}_{\bar{t},\bar{x}}-\Gamma^{\bar{t}}_{\bar{x}x}\Sigma_x^{\bar{t}t}
\bar{\mathcal{J}}_{t,x}\Gamma^{t}_{x\bar{x}}\Sigma_{\bar{x}}^{t\bar{t}})f_z^\tau(\bar{z}), $$
$$f_z^{\mathcal{D}_j\tau}(\bar{z})=\mathcal{D}_jf_z^\tau(\bar{z}),\quad f_z^{\mathcal{I}_0^{ij}(\tau)}(\bar{z})=(\mathcal{I}_0^{ij}+\mathcal{J}_{\bar{t},\bar{x}}^{0,ij}-\Gamma^{\bar{t}}_{\bar{x}x}\Sigma_x^{\bar{t}t}
\bar{\mathcal{J}}_{t,x}^{0,ij}\Gamma^{t}_{x\bar{x}}\Sigma_{\bar{x}}^{t\bar{t}})f_z^\tau(\bar{z}).$$
 It is easy to verify that for $\tau\in \mathcal{F}_F$, $f^\tau_z$ is indeed a $\hat{\mathbf{\Pi}}$-polynomial and for $\tau\in \mathcal{F}_F$
\begin{equation}\hat{\mathcal{Q}}_{<|\tau|_{\mathfrak{s}}}f_z^\tau(z)=0,\forall z\in\mathbb{R}^{d+1},\end{equation}
where $\hat{\mathcal{Q}}_{<\gamma}$ is the projection onto $\hat{{T}}_{<\gamma}$.
Set $\mathcal{T}:=\mathcal{H}_F\cup\{\Xi_j,j=1,2,3\}$.
Define operators $\check{\Gamma}_{x\bar{x}}^t:\mathcal{T}\rightarrow\mathcal{T}, \check{\Sigma}_{x}^{t\bar{t}}:\mathcal{T}\rightarrow\mathcal{T}$ by setting
$$\check{\Gamma}_{x\bar{x}}^t\Xi_j=\Xi_j,\quad \check{\Gamma}_{x\bar{x}}^tX^k=\Gamma_{x\bar{x}}^tX^k,\quad \check{\Sigma}_{x}^{t\bar{t}}\Xi_j=\Xi_j,\quad \check{\Sigma}_{x}^{t\bar{t}}X^k={\Sigma}_{x}^{t\bar{t}}X^k,$$
and then recursively by
$$\check{\Gamma}_{x\bar{x}}^t(\tau\bar{\tau})=(\check{\Gamma}_{x\bar{x}}^t\tau)(\check{\Gamma}_{x\bar{x}}^t\bar{\tau}),\quad \check{\Sigma}_{x}^{t\bar{t}}(\tau\bar{\tau})=(\check{\Sigma}_{x}^{t\bar{t}}\tau)(\check{\Sigma}_{x}^{t\bar{t}}\bar{\tau}),$$
as well as
$$\check{\Gamma}_{x\bar{x}}^t(\mathcal{I}\tau)=\mathcal{I}\check{\Gamma}_{x\bar{x}}^t\tau+(\bar{\mathcal{J}}_{t,x}\Gamma^t_{x\bar{x}}
-\Gamma^t_{x\bar{x}}\bar{\mathcal{J}}_{t,\bar{x}})f^{\tau}_{t,\bar{x}}(t,\bar{x}),$$
$$\check{\Gamma}_{x\bar{x}}^t(\mathcal{I}_0^{ij}\tau)=\mathcal{I}^{ij}_0\check{\Gamma}_{x\bar{x}}^t\tau+(\bar{\mathcal{J}}^{0,ij}_{t,x}\Gamma^t_{x\bar{x}}
-\Gamma^t_{x\bar{x}}\bar{\mathcal{J}}^{0,ij}_{t,\bar{x}})f^{\tau}_{t,\bar{x}}(t,\bar{x}),$$
$$\check{\Sigma}_{x}^{t\bar{t}}(\mathcal{I}\tau)=\mathcal{I}\check{\Sigma}_{x}^{t\bar{t}}\tau+(\bar{\mathcal{J}}_{t,x}\check{\Sigma}_{x}^{t\bar{t}}
-\check{\Sigma}_{x}^{t\bar{t}}\bar{\mathcal{J}}_{\bar{t},x})f^{\tau}_{\bar{t},x}(\bar{t},x),$$
$$\check{\Sigma}_{x}^{t\bar{t}}(\mathcal{I}_0^{ij}\tau)=\mathcal{I}_0^{ij}\check{\Sigma}_{x}^{t\bar{t}}\tau+(\bar{\mathcal{J}}_{t,x}^{0,ij}
\check{\Sigma}_{x}^{t\bar{t}}
-\check{\Sigma}_{x}^{t\bar{t}}\bar{\mathcal{J}}_{\bar{t},x}^{0,ij})f^{\tau}_{\bar{t},x}(\bar{t},x).$$
We also extend the definition of $f_z^\tau$ to all of $\tau\in\mathcal{T}$ by linearity. It is straightforward to verify that we have
\begin{equation}f_{t,x}^{\check{\Gamma}^t_{x\bar{x}}\tau}=f^\tau_{t,\bar{x}},\quad f_{t,x}^{{\check{\Sigma}}_{x}^{t\bar{t}}\tau}=f^\tau_{\bar{t},{x}}. \end{equation}
The map $\mathcal{Z}$ is defined as follows: Given a model $\hat{\mathbf{\Pi}}=(\Pi,\Gamma,\Sigma)\in\hat{\mathcal{M}}_F$, we define a new model $\mathcal{Z}\hat{\mathbf{\Pi}}=(\check{\Pi},\check{\Gamma},\check{\Sigma})$ on $\mathfrak{T}_F$ with
$$\check{\Pi}^t_x\tau=\mathcal{R}_tf_{t,x}^\tau, $$
and $\check{\Gamma},\check{\Sigma}$ being defined as above.
It follows from (3.2) that $\check{\Pi}^t_y=\check{\Pi}^t_x\check{\Gamma}^t_{xy}$. By definition of $\check{\Gamma},\check{\Sigma}$ we can easily verify by recursion that
$\check{\Gamma}^t_{xx}=1 ,\check{\Sigma}_x^{tt}=1$. Since the map $\tau\mapsto f_x^\tau$ is injective, (3.2) also implies that $\Sigma_x^{sr}\Sigma_x^{rt}=\Sigma_x^{st}, \Sigma_x^{st}\Gamma^t_{xy}=\Gamma^s_{xy}\Sigma^{st}_y, \Gamma^t_{xy}\Gamma^t_{yz}=\Gamma^t_{xz}$. It remains to verify that the required analytical bounds hold. The bounds of
 $\check{\Pi}$ follow from the corresponding on $\Pi$ and the fact that $\mathcal{R}_tf^\tau_{t,x}=\Pi_x^tf_{t,x}^\tau(t,x)$ combined with (3.1).
  The bounds of $\check{\Gamma}$ and the bounds of $\check{\Sigma}$ on $\mathcal{I}\tau$  follow inductively from the corresponding bounds on $(\Gamma,\Sigma)$ combined with [14, Lem. 5.21]. The bounds of $\check{\Sigma}$ on $\mathcal{I}_0^{ij}$ follow from the proof in Proposition 2.6 and the bound of $\Sigma$. The fact that the new model is admissible follows from
  $$\check{\Pi}_x^t\mathcal{I}(\tau)=\Pi^t_x(\mathcal{I}+\mathcal{J}_{t,x}-\bar{\mathcal{J}}_{t,x})f^\tau_{t,x}(t,x),$$
   $$\check{\Pi}_x^t\mathcal{I}_0^{ij}(\tau)=\Pi^t_x(\mathcal{I}_0+\mathcal{J}^{0,ij}_{t,x}-\bar{\mathcal{J}}^{0,ij}_{t,x})f^\tau_{t,x}(t,x),$$
  and the definition of $\bar{\mathcal{J}}_{t,x}, \bar{\mathcal{J}}^{0,ij}_{t,x}$.
Setting $\tau(h,\mathbf{\Pi})=\mathcal{Z}(\mathcal{Y}(h,\mathbf{\Pi}))$,
we have
$$\tau(h,\mathfrak{L}(\xi))=\mathfrak{L}(\xi+h),$$
for every smooth $\xi$ and $h$. Moreover,  by similar arguments as in [17, Section 5.1] we have
$$\tau(h,M_g\mathbf{\Pi})=M_g\tau(h,\mathbf{\Pi}),$$
for every smooth $h$ and every smooth $\mathbf{\Pi}\in\mathcal{M}_F$ and every $g\in\mathfrak{R}$.
Now by similar arguments as in the proof of [17, Prop. 4.7],   Assumptions 3 in \cite{HM16} holds.

Furthermore, for any $h\in E$ supported in $(0,1)\times \mathbb{R}^d$ and any model $\mathbf{\Pi}\in\mathcal{M}$, the solution $u^h$ to the first component in (2.7) with model $\mathbf{\Pi}^h=\tau(h,\mathbf{\Pi})$ is related to the solution $\hat{u}^h$ to
$$u^{h,i}=-\frac{1}{2}\sum_{i_1,j=1}^3\mathcal{T}^{ii_1}_j\mathbf{R}^+ (u^{h,i_1}\star u^{h,j})+\sum_{i_1=1}^3\bar{\mathcal{T}}^{ii_1}_j\mathbf{R}^+\Xi_{i_1}+\sum_{i_1=1}^3\mathcal{G}^{ii_1}(u_0^{i_1}+h^{i_1}).$$
with model $\mathbf{\Pi}$ by $\mathcal{R}^hu^h=\mathcal{R}\hat{u}^h$. Here, $\mathcal{R}^h$ and $\mathcal{R}$ denote the reconstruction operators for the model $\mathbf{\Pi}^h$ and $\mathbf{\Pi}$, respectively.  By the implicit function theorem in [7, Theorem 19.28] and similar arguments as in the proof of [3, Proposition 4.7], Assumption 4 follows.

Finally we come to Assumption 5: the map $u_{s,t}$ is Fr\'{e}chet differentiable w.r.t. the initial value. Denote its derivative in the direction $v$ by $J_{s,t}v$. Since the solution to the second component of the solution to the fixed point (2.7) takes values in a sector $V\subset \bar{T}\oplus T_{\geq \zeta}$ with $\zeta=\alpha+3$, $J_{s,t}$ is a bounded linear operator from $\mathcal{C}^\eta$ to $\mathcal{C}^\zeta$. Applying [17, Corollary A.3], the derivative of $u_t$ w.r.t. $h\in L^p([0,1],X_0)$ is given by
$$\mathcal{D}u_t(u_0,\mathbf{\Pi})h=\int_0^t J_{s,t}h(s)ds,$$
for $(\Phi_0,\mathbf{\Pi})\in\mathcal{N}_t$.
By this and the same arguments as in the proof of Theorem 4.8 in \cite{HM16}, the results follow.$\hfill\Box$

\section{Two dimensional case}

We consider (1.1) on $\mathbb{T}^2$.
The strong Feller property in this case can be obtained in the same way as for the three dimensional case. More precisely, we choose $U=(\mathcal{C}^\eta)^2$ for $\eta\in(-\kappa,0)$ for $\kappa$ small enough, and
$$\aligned\mathcal{F}_0=&\{\mathbf{1},  \mathcal{I}^{ii_1}(\Xi_{i_1}), \mathcal{I}^{ii_1}(\Xi_{i_1})\mathcal{I}^{jj_1}(\Xi_{j_1}), i,j,i_1,j_1=1,2,3 \}.\endaligned$$
All the arguments in Section 3 can be applied in this case. It is well-known that the invariant measure of the equation is given by a Gaussian measure $\mu$ (cf. \cite{DD02}). By using the invariant measure in this case we obtain the following result:

\vskip.10in
\th{Theorem 4.1} For the N-S equation driven by space-time white noise on $\mathbb{T}^2$, its solutions are almost surely global in time for every initial condition $u_0\in\mathcal{C}^\eta$ for $\eta\in (-\frac{1}{2},0)$.

\proof We use  $\zeta_{x}$ to denote the explosion time for $u_\cdot(x,\mathbf{\Pi}), x\in \mathcal{C}^\eta$. The invariant measure for the N-S equation driven by space-time white noise on $\mathbb{T}^2$  is given by  Gaussian measure $\mu$ (cf. \cite{DD02}), which has full support. By \cite{DD02} we know that $P(\zeta_x=\infty)=1$ for $\mu$-almost every starting point $x$. By strong Feller property we know that for every $t\geq0$, $x\mapsto P(t<\zeta_x)$ is continuous, which implies that
$P(\zeta_x=\infty)=1$ for every starting point $x$ in $\mathcal{C}^\eta$.

\vskip.10in
\no\textbf{Appendix}

In this section we are  in the setting of Section 3 and we recall the following assumptions from \cite{HM16}. In \cite{HM16}, it has been proved that the Markov semigroup satisfies the strong Feller property under the Assumptions 1-5.
\vskip.10in

\th{Assumption 1} The preimage of $U$ under the map
$$\mathbf{u}: (s,t,u_0,\mathbf{\Pi})\mapsto u_{s,t}(u_0,\mathbf{\Pi}),$$
is open and $\mathbf{u}$ is jointly continuous on $\mathbf{u}^{-1}(U)$. Furthermore, $\mathbf{u}$ is Fr\'{e}chet differentiable in $u_0$ at every point of $\mathbf{u}^{-1}(U)$.
\vskip.10in

Define the sets
$$\mathcal{N}_t=\{(u_0,\mathbf{\Pi}): u_t(u_0,\mathbf{\Pi})\neq \Delta\}.$$
We will denote the Fr\'{e}chet derivative of $u_t$ in the direction $v\in U$ by $Du_t(u_0,\mathbf{\Pi})v$, with the understanding that $Du_t$ is only defined on $\mathcal{N}_t$.
\vskip.10in

\th{Assumption 2} We are given a lower semi-continuous map $r:[0,1]\times U\times \mathcal{M}\rightarrow[0,\infty]$ with the following properties.

1. For every $u_0\in U$ and every $t\in [0,1]$, the map $\omega\mapsto r_t(u_0,\mathbf{Z}(\omega))$ is $\mathcal{F}_t$-measurable.

2. For every $(u_0,\mathbf{\Pi})\in U\times \mathcal{M}$, one has $r_0(u_0,\mathbf{\Pi})<\infty$ and the map $t\mapsto r_t(u_0, \mathbf{\Pi})$ is continuous and increasing.

3. One has $\{(u_0,\mathbf{\Pi}):u_t(u_0,\mathbf{\Pi})=\Delta\}=\{(u_0,\mathbf{\Pi}):r_t(u_0,\mathbf{\Pi})=\infty\}$.

4. For every $t\in (0,1]$, the map $r_t$ is locally Lipschitz continuous on $\mathcal{N}_t$.
\vskip.10in

We consider a space $E=L^p([0,1],X_0)\subset L^2(\mathbb{R},L^2)$, where $X_0$ is some separable Banach subspace of $L^2$ and $p\in (2,\infty)$.
Define $E_s:=L^p([0,s],X_0)$.
\vskip.10in

\th{Assumption 3} We are given a continuous action $\tau: E\times \mathcal{M}\rightarrow\mathcal{M}$ of $E$ onto $\mathcal{M}$ such that, for every $\mathbf{\Pi}\in\mathcal{M}$, the map $h\mapsto \tau(h, \mathbf{\Pi})$ is locally Lipschitz continuous and such that, for every $h\in E$, the identity
$$\mathbf{Z}(\omega+h)=\tau(h,\mathbf{Z}(\omega)),$$
holds $\mathbf{P}$-almost surely. Furthermore, the action $\tau$ is compatible with the maps $u_{s,t}$ in the sense that if $h$ is such that $h(r)=0$ for $r\in [s,t]$, then
$$u_{s,t}(u_0,\tau(h,\mathbf{\Pi}))=u_{s,t}(u_0,\mathbf{\Pi}),$$
for every $u_0\in \bar{U}$.
\vskip.10in

\th{Assumption 4} For every $(u_0,\mathbf{\Pi})\in\mathcal{N}_t$, the map $h\mapsto u_t(u_0,\tau(h,\mathbf{\Pi}))$ is Fr\'{e}chet differentiable at $h=0$.
\vskip.10in

We denote this Fr\'{e}chet derivative by $\mathcal{D}u_t(u_0,\mathbf{\Pi})$. Let $L(U,E_s)$ to denote the space of bounded linear operators $U\rightarrow E_s$.
\vskip.10in

\th{Assumption 5} For every $s\leq t$ with $t\in(0,1]$, we are given a map $A_t^{(s)}:\mathcal{N}_s\rightarrow L(U,E_s)$, and these maps are compatible in the sense that, for any $0<s<r\leq t$, any $(u_0,\mathbf{\Pi})\in\mathcal{N}_r$, and any $v\in U$, one has
$$A_t^{(r)}(u_0,\mathbf{\Pi})v|_{[0,s]}=A_t^{(s)}(u_0,\mathbf{\Pi})v.$$
Furthermore, for every $u_0\in U$ and $s\leq t$, the map $\omega\mapsto A_t^{(s)}(u_0,\mathbf{Z}(\omega))$ is $\mathcal{F}_s$-measurable and one has the identity
$$Du_t(u_0,\mathbf{\Pi})v+\mathcal{D}u_t(u_0,\mathbf{\Pi})( A_t^{(s)}(u_0,\mathbf{\Pi})v)=0,$$
for all $v\in U$ and all $(u_0,\mathbf{\Pi})\in\mathcal{N}_t$. Furthermore, for every $0<s\leq t\leq 1$, the map $A_t^{(s)}$ is locally Lipschitz continuous from $\mathcal{N}_s$ to $L(U,E_s)$ and bounded on every set of the form $\{(u_0,\mathbf{\Pi}):r_s(u_0,\mathbf{\Pi})\leq R\}$ with $R>0$.
\vskip.10in

\end{document}